\documentclass[11pt]{amsart}
\usepackage{verbatim}

\theoremstyle{plain}
\newtheorem{lem}{Lemma}
\newtheorem{thm}[lem]{Theorem}
\newtheorem{prop}[lem]{Proposition}
\newtheorem{cor}[lem]{Corollary}
\theoremstyle{definition}
\newtheorem{rem}[lem]{Remark}
\newtheorem{ex}[lem]{Example}

\begin{document}
\title{Hall-Littlewood vertex operators and
generalized Kostka polynomials}
\author{Mark Shimozono}
\author{Mike Zabrocki}

\maketitle

\pagestyle{plain}

\newcommand{\C}{\mathbb{C}}
\newcommand{\Q}{\mathbb{Q}}
\newcommand{\Z}{\mathbb{Z}}
\newcommand{\N}{\mathbb{N}}
\newcommand{\la}{\lambda}
\newcommand{\La}{\Lambda}
\newcommand{\Roots}{\mathrm{Roots}}
\newcommand{\wt}{\mathrm{wt}}
\newcommand{\inner}[2]{\langle\,#1\,,\,#2\,\rangle}
\newcommand{\w}{v}
\newcommand{\V}{\mathcal V}
\newcommand{\mone}{*}
\newcommand{\cc}{|}
\newcommand{\GGG}{\mathcal{G}}
\newcommand{\FF}{\mathbb{F}}
\newcommand{\HH}{\overline{H}}
\newcommand{\BB}{\overline{B}}

\newcommand{\DIW}[1]{\Z_\ge^{#1}}
\newcommand{\Part}{\mathcal{P}}
\newcommand{\Pa}[1]{{\Part}^{#1}}
\newcommand{\cb}{\overline{c}}
\newcommand{\sign}{\mathrm{sign}}
\newcommand{\coeff}{{\Big{|}}}
\newcommand{\HHH}{\mathbb{H}}
\newcommand{\HHHH}{\mathcal{H}}
\newcommand{\Lah}{\widehat{\La}}
\newcommand{\End}{\mathrm{End}}
\newcommand{\Hom}{\mathrm{Hom}}

\section{Introduction}

Kostka-Folkes polynomials may be considered
as coefficients of the formal power series representing the character
of certain graded $GL(n)$-modules. These $GL(n)$-modules are defined by
twisting the coordinate ring of the nullcone by a suitable line bundle
\cite{Br} and the definition may be generalized by twisting the
coordinate ring of any nilpotent conjugacy closure in $gl(n)$ by a
suitable vector bundle \cite{SW}.  The resulting
polynomials have been called generalized Kostka polynomials.

Jing defined a vertex operator that generates
the Hall-Littlewood symmetric function $Q[X;q]$ \cite{J}, thereby
giving an elegant symmetric function recursion
for the Kostka-Foulkes polynomials. Garsia used a variant
of Jing's vertex operator to derive various new formulas for the
Kostka-Foulkes polynomials.

Our point of departure was the observation that the
Hall-Littlewood vertex operators can be used to obtain formulas for
generalized Kostka polynomials.  Our treatment uses Garsia's plethystic
type formulas.

One striking fact is that the $\Z[q,q^{-1}]$-linear span of $n$-fold
compositions of components of the Hall-Littlewood vertex operators,
is isomorphic to
$K_{G\times \C^*}(\mathcal{N})$, the $GL(n) \times \C^*$-equivariant
$K$-theory of the nullcone.  Under this isomorphism, an $n$-fold
composite operator is sent to the class of
the Euler characteristic of a twisted module.
This fact has a generalization for all the nilpotent conjugacy
class closures in $gl(n)$. 

 These
Grothendieck groups were studied in \cite{KKSW}.
We derive many explicit relations among the vertex operators,
most of which can be in interpreted as relations in the Grothendieck groups
which arise from certain Koszul complexes.  This allows for more
explicit proofs of some basis theorems for these
Grothendieck groups that were proved in \cite{KKSW} using
geometric arguments.

There is a particularly well-behaved subfamily of the
generalized Kostka polynomials, namely, those that
are indexed by a sequence of rectangular partitions.
For this subfamily almost all of the formulas for Kostka-Foulkes
polynomials have generalizations.  There are combinatorial formulas
involving Littlewood-Richardson tableaux \cite{S1} \cite{ScWa},
rigged configurations \cite{KSS}, and inhomogeneous paths with
energy function \cite{S2} \cite{ScWa}.

\section{Hopf algebra of symmetric functions and plethysm}

This section contains standard background material on symmetric
functions which can be found in \cite{M}.  The possible exception
to this is the definition of the plethystic notation used here.

Let $\La=\La_\FF$ be the algebra of symmetric functions over
a field $\FF$ of characteristic zero.  $\La$ may be defined as the
polynomial algebra $\FF[p_1,p_2,\dotsc]$ where the $p_k$ are commuting
algebraically independent variables.  Let $\Part$ be the set of
partitions.  For $\la=(\la_1,\dots,\la_k)\in\Part$, write
$|\la|=\sum_{i=1}^k \la_i$ and let $p_\la = p_{\la_1} \dots p_{\la_k}$
denote the power symmetric function.  Declaring $p_k$ to have degree $k$,
$\La$ is endowed with a grading $\La = \bigoplus_{n\ge 0} \La^n$
where $\La^n$ is the homogeneous component of degree $n$.
$\La^n$ has $\FF$-basis $\{p_\la\mid \la\in\Part, |\la|=n\}$.

$\La$ may be realized by symmetric series.
Let $X=(x_1,x_2,\dotsc)$ be a countable set of commuting indeterminates.
Denote by $\FF[[X]]$ the $\FF$-algebra of formal power series
in the $x_i$ where each $x_i$ has degree 1,
and $\FF_b[[X]]\subset \FF[[X]]$ the subalgebra consisting of series
whose monomials are of bounded degree.
Let $\La^X \subset \FF_b[[X]]$ denote the $\FF$-subalgebra consisting of
the series that are symmetric in the variables $x_i$.  By the universal
mapping property of a polynomial algebra, there is a
$\FF$-algebra homomorphism $\La\rightarrow \La^X$
given by $p_k \mapsto x_1^k+x_2^k+\dotsm$.  This map is
in fact a graded $\FF$-algebra isomorphism.

There is a scalar product on $\La$
for which the power symmetric functions are an orthogonal basis:
\begin{equation*}
  \inner{p_\la}{p_\mu} = \delta_{\la\mu} z_\la
\end{equation*}
where $\delta_{\la\mu}$ is the Kronecker delta,
$z_\la = \prod_i n_i(\la)!\, i^{n_i(\la)}$, and $n_i(\la)$ is the number
of parts of size $i$ in the partition $\la$.

We now define notation for plethystic substitution.

Let $q\in\FF$ be a distinguished element that is transcendental over $\Q$.
Fix an element $E\in\FF_b[[X]]$.  Define $p_k[E]\in \FF_b[[X]]$
to be the series obtained from $E$
by replacing $q$ by $q^k$ and $x_i$ by $x_i^k$ for all $i\ge1$.
By the universal mapping property of a polynomial algebra, there is a
unique $\FF$-algebra homomorphism $\La\rightarrow\FF_b[[X]]$ such that
$p_k \mapsto p_k[E]$.  The image of $P\in\La$ under this map is denoted
$P[E]$.

Switching viewpoints and considering $p_k[E]$ for various $E$, we
see directly from the definitions that
$E\mapsto p_k[E]$ is a $\Q$-algebra homomorphism
$\FF_b[[X]] \rightarrow \FF_b[[X]]$ and $\La^X \rightarrow \La^X$
given by
$p_k[E\pm E'] = p_k[E]\pm p_k[E']$ and $p_k[E E'] = p_k[E]p_k[E']$
for all $E,E'\in\FF_b[[X]]$.
It is \textit{not} an $\FF$-algebra 
homomorphism since it changes the scalar $q$ to $q^k$.

By abuse of notation the variable $X$ will also be used to represent the
formal sum $x_1 + x_2 + x_3 + \dotsm$.  By definition
\begin{equation*}
  p_k[X]=p_k[x_1 + x_2 + x_3 + \dotsm] = x_1^k + x_2^k + x_3^k +\dotsm.
\end{equation*}

\begin{ex} $e_2 = p_1^2/2 - p_2/2 \in \Lambda$.  For the expression
$E = x_1 + x_2$ we have
$$e_2[x_1+x_2] = (x_1+x_2)^2/2 - (x_1^2 + x_2^2)/2 = x_1 x_2$$
If $E = \frac{1}{1-q}$ then $p_1[\frac{1}{1-q}]=\frac{1}{1-q}$ and
$p_2[\frac{1}{1-q}]=\frac{1}{1-q^2}$, so that
$$e_2 \left[ \frac{1}{1-q} \right] = \frac{1}{2(1-q)^2} -
\frac{1}{2(1-q^2)}= \frac{q}{(1-q)(1-q^2)}$$
\end{ex}

Recall that each $E\in\FF_b[[X]]$ defines an $\FF$-algebra homomorphism
$\La\rightarrow\FF_b[[X]]$ such that $P\mapsto P[E]$.
\begin{enumerate}
\item If $E=X$ then $P\mapsto P[X]$ yields the above isomorphism
$\La \cong \La^X$.
\item If $E=0$ then the $\FF$-algebra homomorphism
$\epsilon:\La\rightarrow\FF$ given by $P\mapsto P[0]$
is the \textit{counit} for $\La$; it selects the constant term.
\item If $E=x_1+x_2+\dots+x_m$ then $P\mapsto P[x_1+x_2+\dots+x_m]$
is the $\FF$-algebra epimorphism $\La\rightarrow \FF[x_1,\dots,x_m]^{S_m}$
onto the $\FF$-algebra of symmetric polynomials in $x_1$ through $x_m$
where $S_m$ is the symmetric group on $m$ letters.
\item If $E = X/(q-1)$ then the map $P[X]\mapsto P[X/(q-1)]$ is an
$\FF$-algebra isomorphism $\La\rightarrow\La$ with inverse map
$P[X]\mapsto P[X(q-1)]$.  In particular if $\{f_\la\}$ is a basis of
$\La$ then so are $\{f_\la[X(q-1)]\}$ and $\{f_\la[X/(q-1)]\}$.
\end{enumerate}

We now discuss the Hopf algebra structure on $\La$.
Let $Y=y_1+y_2+\dotsm$ where the $y_i$ are another countable
set of commuting indeterminates.  Let $\La^{X;Y}$ denote the
$\FF$-subalgebra of $\FF_b[[X,Y]]$ consisting of the series that
are symmetric in both the $x_i$ and the $y_j$.  Then there is an
isomorphism $\La \otimes_\FF \La \rightarrow \La^{X;Y}$ given by
$f \otimes g \mapsto f[X] g[Y]$ for $f,g\in\La$.  Letting $E=X+Y$,
the map $P\mapsto P[X+Y]$ defines an $\FF$-algebra
homomorphism $\Delta:\La \rightarrow \La \otimes \La$.
$\Delta$ is the comultiplication map for $\La$.  
The maps $\Delta$ and $\epsilon$ give $\La$ the structure of a
coassociative cocommutative Hopf algebra.

There is a scalar product on $\La \otimes \La$ defined by
\begin{equation*}
  \inner{f_1[X] g_1[Y]}{f_2[X] g_2[Y]}_{\La^X\otimes \La^Y} =
  \inner{f_1}{f_2} \inner{g_1}{g_2}
\end{equation*}
for all $f_i,g_i\in \La$ for $i=1,2$ and this relation is then extended by
linearity. With respect to this scalar product on $\La \otimes \La$,
the map $\Delta$ is adjoint to the multiplication map
$\La \otimes \La \rightarrow \La$ given by $f \otimes g \mapsto f g$.
That is,
\begin{equation} \label{eq:scalar coprod}
  \inner{f}{g h} = \inner{f[X+Y]}{g[X] h[Y]}
\end{equation}
for all $f,g,h\in\Lambda$.

To deal with the Cauchy kernel it is necessary to work in a
completion of $\La$.  Let $\Lah = \FF[[p_1,p_2,\dotsc]]$,
the $\FF$-algebra of formal power series in the $p_i$.
The symmetric realization $\Lah^X$ of $\Lah$ is given by
the $\FF$-subalgebra of symmetric series
in $\FF[[X]]$.  Given any element $E\in\FF_b[[X]]$ 
with no constant term and $P\in\Lah$, the plethysm
$P[E]\in\FF[[X]]$ may be defined by substitution as before.

Define the Cauchy kernel $\Omega\in\Lah$ by
\begin{equation*}
  \Omega =  \exp\left(\sum_{r\ge1} p_r/r\right).
\end{equation*}
The following formulas are consequences of the definitions:
\begin{equation*}
\begin{split}
\Omega[ zX ] &= \prod_i \frac{1}{1-z x_i} = \sum_{k \geq 0} z^k h_k[X] \\
\Omega[-zX] &= \prod_i 1-z x_i = \sum_{k \geq 0} (-z)^k e_k[X] \\
\Omega[X+Y] &= \Omega[X] \Omega[Y] \\
\Omega[X-Y] &= \Omega[X] \Omega[-Y] = \Omega[X]/\Omega[Y] \\
\Omega[XY] &= \sum_{\la\in\Part} s_\la[X] s_\la[Y]
\end{split}
\end{equation*}
where $h_k$, $e_k$, and $s_\la$ are the homogeneous, elementary,
and Schur symmetric functions respectively.

Two bases $\{a_\la \}$ and $\{b_\la\}$ of $\La$ are dual
with respect to the scalar product if and only if
$\sum_{\la\in\Part} a_\la[X] b_\la[Y] = \Omega[XY]$.

For $f\in\La$, let $f^\perp$ be the linear operator on $\La$
that is adjoint to multiplication by $f$ with respect to the scalar product:
\begin{equation} \label{eq:perp def}
  \inner{f^\perp(g)}{h} = \inner{g}{f h}
\end{equation}
for all $g,h \in \La$.  This operator is usually referred to as ``f-perp''
or ``skewing by f.''

The skewing operators and the map $\Delta$ (or equivalently the
plethystic map $P\mapsto P[X+Y]$) can be expressed in terms of each other.
For any $f,g,h\in \La$,
\begin{equation} \label{eq:skew pleth}
  \inner{f^\perp(g)}{h} = \inner{g[X+Y]}{f[Y] h[X]}
\end{equation}
using \eqref{eq:perp def}, \eqref{eq:scalar coprod},
and the commutativity of multiplication in $\La$.
Let $\{a_\la\}$ and $\{b_\la\}$ be dual bases of $\La$.
Then for every $P\in\La$,
\begin{equation} \label{eq:coprod perp}
  P[X+Y] = \sum_\la (b_\la^\perp(P))[X] a_\la[Y].
\end{equation}

\section{Definition of the operator}

Using plethystic substitution we define a family of linear operators on
symmetric functions.  Define the formal Laurent series
$H(Z^k)$ in an ordered set of variables $Z^k=(z_1,z_2,\dotsc,z_k)$
with coefficients given by operators on $\La$, 
which acts on $P \in \La$ by
\begin{equation} \label{eq:H gen def}
H(Z^k) P[X] = P[X - (1-q)Z^\mone] \Omega[ ZX ] R(Z^k)
\end{equation}
where $Z^\mone = \sum_{i=1}^k z_i^{-1}$, $Z = \sum_{i=1}^k
z_i$ and
$R(Z^k) = \prod_{1 \leq i < j \leq k} 1- z_j/z_i$.
For $\w \in \Z^k$, define the operator $H^q_\w$ by
\begin{equation} \label{eq:H def}
  H^q_{\w} P[X] = H(Z^k) P[X]\coeff_{z^{\w}}.
\end{equation}

\begin{rem}
\begin{enumerate}
\item If $k=1$ this is Garsia's \cite{G1} \cite{G2} version
of Jing's Hall-Littlewood vertex operator \cite{J}.
\item At $q=0$ this formula reduces to plethystic notation for a repeated
application of the Schur function vertex operator that is due to Bernstein
\cite[p.96, \# 29 (d)]{M}, \cite{G2}.  So 
$H_{(\w_1)}^0 H_{(\w_2)}^0 \cdots H_{(\w_k)}^0 = H_{\w}^0$
for any $\w\in\Z^k$.
If $\w$ is a partition $\la=(\la_1,\dots,\la_k)$ with at most $k$ parts
then $H_{\la}^0 \, 1 = s_\la[X]$.
\item When $q=1$ and $\la$ is a partition, this formula reduces to
multiplication by the Schur function $s_\la$:
$\Omega[ZX]R(Z^k) \coeff_{z^\la} = s_{\la}[X].$
\item $H^q_\la\, 1 = s_\la$ for any partition $\la=(\la_1,\dots,\la_k)$
with at most $k$ parts.  If $\w_k<0$ then $H^q_\w \, 1 = 0$.
\end{enumerate}
\end{rem}

This operator possesses the same shifted skew symmetry in its
subscript that Schur functions have.
Let $\DIW{k}$ denote the set of dominant integral weights in $\Z^k$,
that is, the weakly decreasing sequences in $\Z^k$.
Let $\Pa{k}\subset\DIW{k}$ be the set of partitions of length at most $k$,
which are always regarded as having $k$ parts, some of which may be zero.

\begin{prop} Let $\w \in \Z^k$.  Then
\begin{equation*}
H^q_\w = - H^q_{(\w_1,\w_2 \ldots,\w_{i+1}-1, \w_i+1, \dotsc, \w_k)}.
\end{equation*}
In particular if $\w_{i+1} = \w_i + 1$ then $H^q_\w = 0$.
\end{prop}
\begin{proof} Let $s_i Z^k$ be the sequence $Z^k$ with
$z_i$ and $z_{i+1}$ exchanged.  Then
$H(s_i Z^k) = 
\frac{1-z_{i}/z_{i+1}}{1-z_{i+1}/z_{i}} H(Z^k) =
- \frac{z_{i}}{z_{i+1}} H(Z^k)$.  The result follows by taking
the coefficient of $z^\w$ on both sides of this equation.
If $\w_{i+1} = \w_i + 1$ then $H^q_\w = - H^q_\w$ and hence must be zero.
\end{proof}

The following corollary shows that for every $\w\in\Z^k$,
$H^q_\w$ is either zero or up to sign equal to $H^q_\nu$ for some
$\nu\in\DIW{k}$.

\begin{cor} \label{co:Bott} Let $\w \in \Z^k$, $\sigma\in S_k$, and
$\rho=\rho^{(k)} = (k-1, k-2, \ldots, 0)$.  Then
\begin{equation}
  H^q_\w = sign(\sigma) H^q_{\sigma(\w + \rho) - \rho}.
\end{equation}
In particular, if $\w+\rho$ has two equal parts then $H^q_\w = 0$.
Otherwise there is a unique $\sigma\in S_k$ such that
$\sigma(\w+\rho)$ is strictly decreasing, so that
$\sigma(\w+\rho)-\rho\in\DIW{k}$.
\end{cor}

\section{Another formula for the vertex operator}

We derive another formula for the operator $H^q_\nu$, which
will be used to prove a strong linear independence property of the
operators.  To do this, it is convenient to introduce some notation
for the irreducible rational characters of $GL(k)$.  The polynomial
representation ring of $GL(k)$ is isomorphic to the $\FF$-algebra
$\FF[z_1,z_2,\dots,z_k]^{S_k}$ of symmetric polynomials in the
variables $z_i$.  The rational representation
ring $R(GL(k))$ is isomorphic to the $\FF$-algebra of symmetric Laurent
polynomials in the variables $z_1,z_2,\dots,z_k$ or equivalently
the localization $\FF[z_1,\dots,z_k]^{S_k}[\det[Z]^{-1}]$
of the algebra of symmetric polynomials obtained by inverting the
character $\det[Z] = z_1 z_2 \dotsm z_k$ of the determinant module.
$R(GL(k))$ has a distinguished basis $\{ s_\la[Z] \mid \la\in\DIW{k}\}$
where $s_\la[Z]$ is the character of the irreducible finite-dimensional
$GL(n)$-module $V^\la$ of highest weight $\la$.
Explicitly, if $m$ is an integer such that $m \ge -\la_k$ then
$s_\la[Z] = \det[Z]^{-m} s_{\la+(m^k)}[Z]$ where
$m^k=(m,m,\dots,m)\in\Z^k$, $\la+(m^k)$ is the vector sum in $\Z^k$
(which is a partition) and
$s_{\la+(m^k)}[Z]$ is the Schur polynomial.
For $\w\in\Z^k$, define its dual weight
$\w^* = (-\w_k,-\w_{k-1},\dotsc,-\w_1)$.  If $\la\in\DIW{k}$
then so is $\la^*$; it satisfies $(V^\la)^* \cong V^{\la^*}$
where $V^*$ is the contragredient dual of $V$.

Denote by $\cb^\la_{\mu\nu}$ the tensor product multiplicities:
\begin{equation*}
  s_\mu[Z] s_\nu[Z] = \sum_{\la\in\DIW{k}} \cb^\la_{\mu\nu} s_\la[Z]
\end{equation*}
for $\la,\mu,\nu\in\DIW{k}$.  It is well-known that
\begin{equation} \label{eq:triple mults}
  \cb^\la_{\mu\nu} = \dim (V^{\la^*} \otimes V^\mu \otimes V^\nu)^{GL(n)}
\end{equation}
where $V^G$ denotes the submodule of $V$ fixed by $G$.  Moreover for any
finite-dimensional $G$-module,
\begin{equation} \label{eq:dual mults}
  \dim((V^*)^G) = \dim(V^G).
\end{equation}

Let
\begin{equation} \label{eq:H gen again}
  \HH(Z) P[X] = P[X-(1-q)Z^\mone] \Omega[ZX].
\end{equation}
Then by \eqref{eq:H gen def} and \eqref{eq:H def},
\begin{equation} \label{eq:H again}
  \HH(Z) = \sum_{\nu\in\DIW{k}} s_\nu[Z] H^q_\nu.
\end{equation}

\begin{prop} \label{pp:old style} For $\nu\in\DIW{k}$,
\begin{equation} \label{eq: old style}
  H^q_\nu = \sum_{\la,\mu\in\Pa{k}} \cb^\la_{\mu\nu}
    s_\la[X] s_\mu[X(q-1)]^\perp.
\end{equation}
\end{prop}
\begin{proof} Let $Y$ be another set of variables and let $\HH(Z)$ act
on the $X$ variables.  We have
\begin{equation*}
\begin{split}
\HH(Z) \Omega[X Y] &= \Omega[X Y] \Omega[(q-1)Z^* Y] \Omega[X Z] \\
&= \Omega[X Y]  \sum_{\nu\in\DIW{k}} \sum_{\la,\mu\in\Pa{k}} s_\nu[Z]
\cb^\nu_{\mu^*\la} s_\mu[(q-1)Y] s_\la[X].
\end{split}
\end{equation*}
The tensor product multiplicity can be rewritten using
\eqref{eq:triple mults} and \eqref{eq:dual mults}
 so that $\cb^\nu_{\mu^*\la} = \cb^\la_{\mu\nu}$.
Clearly we have for fixed $\gamma \in \Pa{k}$ and $\nu \in \DIW{k}$
\begin{align*}
H^q_\nu s_\gamma[X] &= \HH(Z) \Omega[X Y]\coeff_{s_\nu[Z]s_\gamma[Y]} \\
&=  \sum_{\la,\mu\in\Pa{k}}
\cb^\la_{\mu\nu} \left< s_\mu[(q-1)Y]^\perp s_\gamma[Y],\Omega[X Y]
s_\la[X] \right>_Y \\
&=  \sum_{\la,\mu\in\Pa{k}}
\cb^\la_{\mu\nu} s_\la[X] s_\mu[(q-1)X]^\perp s_\gamma[X]
\end{align*}

It is seen that \eqref{eq: old style} holds
on the Schur basis and hence on $\La$.
\end{proof}

\section{Linear independence}

The full strength of the following result is required later,
when we work with operators on $\La$ that are infinite linear combinations
of the $H^q_\nu$.

\begin{prop} \label{pp:indp} Let $\{c_\nu\in\FF\mid \nu\in\DIW{k}, |\nu|
= d \}$ be an arbitrary collection of scalars.  Then the map
$F = \sum_{\nu\in\DIW{k}, |\nu| = d} c_\nu H^q_\nu$ is a well-defined
linear operator on $\La$.  Moreover this operator is $0$ if and only
if $c_\nu=0$ for all $\nu\in\DIW{k}$ with $|\nu| =d$.
\end{prop}
\begin{proof} 

Any infinite linear combination of operators of the
form $s_\la[X] s_\mu[X(q-1)]^\perp$ with $|\la| = |\mu| + d$ is
well-defined because its action on any given element of the basis
$\{s_\la[X/(q-1)]\}$ has only finitely many nonzero summands.

The operator $F = \sum_{\nu\in\DIW{k}, |\nu| = d} c_\nu H^q_\nu$
is well-defined because for fixed partitions $\la$ and $\mu$,
$s_\la[X] s_\mu[X(q-1)]^\perp$ appears in the formula
given in Proposition \ref{pp:old style} for only finitely many
$H^q_\nu$, namely for those such that $\cb^{\la}_{\mu\nu}\not=0$.

For $\gamma\in\DIW{k}$, set $\alpha(\gamma)\in\Pa{k}$ to be the partition
with $\alpha(\gamma)_i = \max(\gamma_i,0)$ and
$\beta(\gamma)\in\Pa{k}$ be the partition defined by $\beta(\gamma)_i = -
min(\gamma_{k+1 - i},0)$.

By Proposition \ref{pp:old style}, for $\tau\in\Pa{k}$,
\begin{equation}
 H^q_{\gamma}(s_{\tau}[X/(q-1)]) =\sum_{\la,\mu\in\Pa{k}} \cb^\la_{\mu\gamma}
    s_\la[X] s_\mu[X(q-1)]^\perp (s_{\tau}[X/(q-1)]).
\end{equation}
The Littlewood Richardson rule implies that the coefficient 
$\cb^\la_{\mu\gamma}=0$ unless $\mu \supseteq \beta(\gamma)$.
Since $\{s_\mu[X(q-1)]\}$ and $\{s_\tau[X/(q-1)]\}$ are dual bases,
we calculate that $H^q_{\gamma}(s_{\tau}[X/(q-1)])=0$ if 
$|\tau| \leq |\beta(\gamma)|$
and $\tau \neq \beta(\gamma)$, and 
$H^q_\gamma (s_{\beta(\gamma)}[X/(q-1)])
 = s_{\alpha(\gamma)}[X]$.  

It follows now that $F=0$ if and only if $c_\gamma=0$ for all $\gamma$.
For if all of the coefficients are not zero, then $\gamma$ is chosen such
that $c_\gamma$ is non-zero and $|\beta(\gamma)|$ is a minimum.  We see then
that 
\begin{equation}
F(s_{\beta(\gamma)}[X/(q-1)]) = \sum_{\beta(\nu) = \beta(\gamma)}
c_\nu s_{\alpha(\nu)}[X] \not= 0
\end{equation}
\end{proof}

\section{Connection with generalized Kostka polynomials}

Next it is shown that the operators $H^q_\nu$ have the same
relation to the generalized Kostka polynomials of \cite{SW} that
the components of Garsia's modified Hall-Littlewood vertex operator
has to the Kostka-Foulkes polynomials.

Let us recall the definition of the generalized Kostka polynomials.
Let $\eta = (\eta_1, \eta_2, \dotsc, \eta_t)$ be a
sequence of positive integers summing to $n$.  Let
\begin{equation*}
  \Roots_\eta = \{ (i,j) \mid 1\le i\le \eta_1+\dots+\eta_k
        < j \le n \text{ for some $k$ }\}
\end{equation*}
be the set of strictly upper block triangular positions
in an $n\times n$ matrix with diagonal block sizes given by $\eta$.
Let $Z$ be the sequence of variables $z_1$ through $z_n$.
Define the formal power series $B_\eta[Z;q]$ by
\begin{equation*}
  B_\eta[Z;q] = \prod_{(i,j)\in\Roots_\eta} \dfrac{1}{1-q z_i/z_j}.
\end{equation*}
Let $J$ be the antisymmetrizer
$J = \sum_{\sigma\in S_n} \sign(\sigma) \sigma$.
Define the linear operator
$\pi:\FF[z_1^\pm,\dots,z_n^\pm]\rightarrow R(GL(n))$ given by
$\pi(f) = J(z^{\rho})^{-1} J(z^\rho f)$.

Let $\gamma \in \Z^n$.  Define the generating series
\begin{equation*}
  \HHH_{\eta,\gamma}[Z;q] = \pi \left( z^\gamma B_\eta[Z;q] \right).
\end{equation*}
Since $\{s_\la[Z]\mid \la\in\DIW{n}\}$ is a basis of $R(GL(n))$,
define $K_{\la\gamma\eta}(q) \in \Z[[q]]$ by
\begin{equation*}
  \HHH_{\eta,\gamma}[Z;q] = \sum_{\la\in\DIW{n}}
  s_\la[Z] K_{\la\gamma\eta}(q).
\end{equation*}
It is shown in \cite{SW} that $K_{\la\gamma\eta}(q) \in \Z[q]$.

By developing the product of geometric series in $B_\eta[Z;q]$
and using the shifted skew symmetry of the index for irreducible
characters
\begin{equation*}
  s_{\w}[Z] = \sign(\sigma) s_{\sigma(\w+\rho)-\rho}[Z]
\end{equation*}
for $\w\in\Z^n$ and $\sigma\in S_n$,
one obtains the following expression \cite{SW}:
\begin{equation} \label{eq:Kostant}
  K_{\la\gamma\eta}(q) =
  \sum_{\sigma\in S_n} \sign(\sigma)
  \sum_{\substack{m:\Roots_\eta\rightarrow\N \\
  \wt(m)=\sigma(\la+\rho)-(\gamma+\rho)}} q^{|m|}
\end{equation}
where $\epsilon_1$ through $\epsilon_n$ is the standard basis of $\Z^n$,
$|m| = \sum_{(i,j)\in\Roots_\eta} m(i,j)$, and
$\wt(m) = \sum_{(i,j)\in\Roots_\eta} m(i,j) (\epsilon_i-\epsilon_j)$.

Given $\eta$ and $\gamma$, write
$\gamma = (\gamma^{(1)}, \gamma^{(2)}, \dotsc, \gamma^{(t)})$
with $\gamma^{(i)}\in\Z^{\eta_i}$ and assume that
$\gamma^{(i)}\in\DIW{\eta_i}$.

\begin{prop} \label{pp: gen coefs}
\begin{equation}\label{eq:genfun}
H^q_{\gamma^{(1)}} H^q_{\gamma^{(2)}} \cdots H^q_{\gamma^{(t)}} =
 \sum_{\la\in\DIW{n}} K_{\la,\gamma,\eta}(q) H^q_\la.
\end{equation}
\end{prop}
\begin{proof} Let $Z^n$ be the entire set of $n$ variables $z_1$
through $z_n$.  Break this set into $t$ collections of successive variables
such that the $i$-th collection
$Z^{(i)}$ has size $\eta_i$ for $1\le i\le t$.  We have the relation
$H(U^k) H(V^\ell) =  \Omega[q U^{\mone} V] H(U^k, V^\ell),$ which
is verified by showing it holds when evaluated on an arbitrary
$P\in\La$:
\begin{align*}
H(U^k) H(V^\ell) P[X] &= P[X - (1-q)(U^{\mone}+V^{\mone})]
\Omega[UX] \\&\hskip .2in\Omega[V(X - (1-q)U^{\mone})] R(U^k) R(V^\ell) \\
&= P[X - (1-q)(U^{\mone}+V^{\mone})]
\Omega[(U+V)X] \\&\hskip .2in\Omega[q U^{\mone}V] R(U^k, V^\ell)
\end{align*}
It follows that
\begin{equation} \label{eq:comp}
H(Z^{(1)}) H(Z^{(2)}) \cdots H(Z^{(t)}) 
= H(Z^n) \prod_{1\leq i<j\leq t} \Omega[q (Z^{(i)})^{\mone} Z^{(j)}]
\end{equation}
Observe that
$H^q_{\gamma^{(1)}} H^q_{\gamma^{(2)}} \cdots H^q_{\gamma^{(t)}}$
is the coefficient of $z^\gamma$ on the left hand side of
\eqref{eq:comp}.  The coefficient of $z^\gamma$ on the right hand side
is $\sum_m q^{|m|} H^q_{\gamma+\wt(m)}$ where $m$ runs over the functions
$m:\Roots_\eta\rightarrow\N$.  By Proposition \ref{pp:indp}
the operators $H^q_\la$ are independent for $\la\in\DIW{n}$.
Hence we may take the coefficient of $H^q_\la$ on both sides.
Using Corollary \ref{co:Bott} the right hand side becomes precisely the
the expression \eqref{eq:Kostant} of $K_{\la\gamma\eta}(q)$.
\end{proof}

If $\gamma$ is such that $\gamma^{(i)}\in\Pa{\eta_i}$, define the
symmetric function
\begin{equation} \label{eq:fin H}
  \HHHH_{(\gamma^{(1)},\dotsc,\gamma^{(t)})}[X;q] =
   H^q_{\gamma^{(1)}} H^q_{\gamma^{(2)}} \cdots H^q_{\gamma^{(t)}} \,1.
\end{equation}
Then
\begin{equation} \label{eq:fin K}
  \HHHH_{(\gamma^{(1)},\dotsc,\gamma^{(t)})}[X;q] =
  \sum_{\la\in\Pa{n}} K_{\la\gamma\eta}(q) s_\la[X]
\end{equation}
which is a finite sum, unlike the expansion \eqref{eq:genfun}
which is an infinite sum.  At $q=1$ this is the
expansion of the product of Schur functions
$s_{\gamma^{(1)}} \dotsm s_{\gamma^{(t)}}$ in the Schur basis.

We have shown that the coefficients of the expansions of
$\HHH_{\eta,\gamma}[Z;q]$ in terms of irreducible characters,
and $H^q_{\gamma^{(1)}} \dots H^q_{\gamma^{(t)}}$ in terms
of the $H^q_\la$ for $\la\in\DIW{n}$, are the same.
Consequently, at least on the level of characters, the questions of
\cite{KKSW} regarding a certain Grothendieck ring of graded modules,
can be translated into questions regarding the span of the
above operators.

\section{Spaces of vertex operators and commutation relations}

By manipulating the order of the variables we may derive several sorts of
explicit commutation relations.  These are applied to prove basis theorems
for spaces of our operators.  Let $\GGG(k,n)$ denote the
$\Z[q]$-span of operators of the form $H^q_\mu H^q_\nu$ where
$\mu\in\Z^k$ and $\nu\in\Z^n$, and $\GGG_\ge(k,n)$ the $\Z[q]$-span
of such operators with the concatenated weight dominant, that is,
$(\mu,\nu)\in\DIW{k+n}$.

The following are vertex operator analogues
of results proven in \cite{KKSW} for Grothendieck groups.  Our proofs
have the advantage of being explicit and working over $\Z[q]$ instead
of $\Z[q,q^{-1}]$.  Moreover our relations among vertex operators
can be lifted to relations in the Grothendieck groups.

\begin{thm} \label{th:dom} $\GGG(k,n)=\GGG_\ge(k,n)$.
\end{thm}

\begin{thm} \label{th:support} If $k>n$ then
$\GGG(k,n)\subset \GGG(k-1,n+1)$ and
$\GGG(n,k)\subset \GGG(n+1,k-1)$.
\end{thm}

\begin{thm} \label{th:comm} $\GGG(k,n)=\GGG(n,k)$.
\end{thm}

Consider the following generating
function equation that follows easily from manipulating the
operators on an arbitrary symmetric function.

\begin{equation} \label{eq:comm}
\begin{split}
H(&U^k,V^\ell) H(W^m Z^n) \Omega[ -q(U^{\mone}W + V^{\mone}W+V^{\mone} Z)] \\
&= 
H(U^k,W^m) H(V^\ell,Z^n) \Omega[-q(U^{\mone} V + W^{\mone}V+W^{\mone}Z)]
(-1)^{\ell m} \prod_{i=1}^m w_i^{\ell} \prod_{i=1}^\ell v_i^{-m}
\end{split}
\end{equation}

Setting $\ell=m=1$ in this formula gives enough relations to prove
Theorem \ref{th:dom}.  For brevity of notation,
define $\cc^\alpha_\beta = |\alpha|-|\beta|$.
If $\mu,\nu\in\DIW{n}$
then we will say that $\nu \slash \mu$ is a vertical strip,
and denote this by $\nu \slash \mu \in \V$, provided that
$\nu_i -\mu_i \in\{ 0,1 \}$ for all $1\le i\le n$.

\begin{lem} For all $a,b\in\Z$ and 
$\mu\in\DIW{k}$ and $\nu\in\DIW{n}$,
\begin{align}
&\sum_{\alpha \slash \mu \in \V} \sum_{\nu \slash \beta \in \V}
(-q)^{\cc^\alpha_\mu+\cc^\nu_\beta} (H_{(\alpha,a+\cc^\nu_\beta)}^q
H_{(b-\cc^\alpha_\mu,\beta)}^q - q H_{(\alpha,a+\cc^\nu_\beta+1)}^q
H_{(b-\cc^\alpha_\mu-1, \beta)}^q )= \nonumber\\
&~\sum_{\alpha \slash \mu \in \V} \sum_{\nu \slash \beta \in \V}
(-q)^{\cc^\alpha_\mu+\cc^\nu_\beta} (q H_{(\alpha,b+\cc^\nu_\beta)}^q
H_{(a-\cc^\alpha_\mu,\beta)}^q - H_{(\alpha,b+\cc^\nu_\beta-1)}^q
H_{(a-\cc^\alpha_\mu+1, \beta)}^q ) \label{eq:com1}
\end{align}
\end{lem}
\begin{proof} From equation \eqref{eq:comm},
\begin{equation*}
\begin{split}
&H(U^k,v) H(w, Z^n) \Omega[-q(U^{\mone}w + v^{\mone}w+v^{\mone} Z)] \\
=\,&H(U^k,w) H(v,Z^n) \Omega[-q(U^{\mone} v + w^{\mone}v+w^{\mone}Z)] (-w/v)
\end{split}
\end{equation*}
The desired identity is obtained by
taking the coefficient of $u^\mu v^a w^b z^\nu$.
\end{proof}

Note that if $b=a+1$, then \eqref{eq:com1} reduces significantly to the
following identity.
\begin{align}
&\sum_{\alpha \slash \mu \in \V} \sum_{\nu \slash \beta \in \V}
(-q)^{\cc^\alpha_\mu+\cc^\nu_\beta} H_{(\alpha,a+\cc^\nu_\beta)}^q
H_{(a+1-\cc^\alpha_\mu,\beta)}^q \nonumber\\
&= \sum_{\alpha \slash \mu \in \V} \sum_{\nu \slash \beta \in \V}
q (-q)^{\cc^\alpha_\mu+\cc^\nu_\beta}  H_{(\alpha,a+\cc^\nu_\beta+1)}^q
H_{(a-\cc^\alpha_\mu, \beta)}^q  \label{eq:com2}
\end{align}

Relations \eqref{eq:com1} and \eqref{eq:com2} suffice to prove
Theorem \ref{th:dom}.  Before giving the proof an example is helpful.

\begin{ex} We wish to show that $H_{(22)}^q H_{(41)}^q$ is a
linear combination of $H_\alpha^q H_\beta^q$ with $\alpha,\beta \in \Z^2$
such that $(\alpha,\beta)$ is dominant.
\begin{align*}
&H_{(22)}^q H_{(41)}^q - q H_{(23)}^q H_{(40)}^q - q
H_{(32)}^q H_{(31)}^q + q^2 H_{(33)}^q H_{(30)}^q \\
&-q H_{(23)}^q H_{(31)}^q + q^2 H_{(24)}^q H_{(30)}^q + q^2
H_{(33)}^q H_{(21)}^q - q^3 H_{(34)}^q H_{(20)}^q = \\
&\hskip .2in q H_{(24)}^q H_{(21)}^q - q^2 H_{(25)}^q H_{(20)}^q
- q^2 H_{(34)}^q H_{(11)}^q + q^3 H_{(35)}^q H_{(10)}^q \\
&\hskip .2in - H_{(23)}^q H_{(31)}^q + q H_{(24)}^q H_{(30)}^q
+ q H_{(33)}^q H_{(21)}^q - q^2 H_{(34)}^q H_{(20)}^q
\end{align*}

By Corollary \ref{co:Bott} many terms vanish and others cancel
with each other.
The terms $H_{(23)}^q H_{(40)}^q$, $H_{(23)}^q
H_{(31)}^q$, $H_{(34)}^q H_{(20)}^q$, and $H_{(34)}^q H_{(11)}^q$
are all zero.  The terms $q^2 H_{(24)}^q H_{(30)}^q$ and $q^2
H_{(33)}^q H_{(30)}^q$ cancel and so do
$q H_{(24)}^q H_{(21)}^q$ and $q H_{(33)}^q H_{(21)}^q$.  When
this relation is reduced and
$H_{(22)}^q H_{(41)}^q$ is expressed alone
on the left hand side of the equation we have
\begin{align*}
H_{(22)}^q H_{(41)}^q = \,&q H_{(32)}^q H_{(31)}^q  
  - q^2 H_{(33)}^q H_{(21)}^q  + q^2 H_{(43)}^q H_{(20)}^q\\
& - q^3 H_{(44)}^q H_{(10)}^q - q H_{(33)}^q H_{(30)}^q 
\end{align*}
On the right hand side of this equation, only $H_{(32)}^q H_{(31)}^q$
does not have the property that the concatenated indexing weights
are dominant.  We apply relation \eqref{eq:com2} to this operator:
\begin{align*}
&H_{(32)}^q H_{(31)}^q - q H_{(33)}^q H_{(30)}^q - q H_{(42)}^q H_{(21)}^q
+ q^2 H_{(43)}^q H_{(20)}^q =\\
&\hskip .2in q H_{(33)}^q H_{(21)}^q - q^2 H_{(34)}^q H_{(20)}^q - q^2
H_{(43)}^q H_{(11)}^q + q^3 H_{(44)}^q H_{(10)}^q
\end{align*}
Therefore
\begin{align*}
H_{(22)}^q H_{(41)}^q = \,&(q^2-q) H_{(33)}^q H_{(30)}^q 
+ q^2 H_{(42)}^q H_{(21)}^q - q^3 H_{(43)}^q H_{(11)}^q \\
&+ (q^2- q^3) H_{(43)}^q H_{(20)}^q + (q^4-q^3) H_{(44)}^q H_{(10)}^q    
\end{align*}
\end{ex}

We now give the proof of Theorem \ref{th:dom}.
\begin{proof} It is enough to show that
if $(\mu,a)\in\DIW{k}$ and $(b,\nu)\in\DIW{n}$ then
\begin{equation*}
H_{(\mu,a)}^q H_{(b,\nu)}^q \in \GGG_\ge(k,n).
\end{equation*}
If $b=a+1$ then in the relation \eqref{eq:com2} the only term
that is not indexed by weights such that
$(\alpha,\beta)\in\DIW{k+n}$, is
$H_{(\mu,a)}^q H_{(a+1,\nu)}^q$, and hence it is in $\GGG_\ge(k,n)$.

Otherwise assume that $b>a+1$.  All of the terms of the equation
\eqref{eq:com1} are of the form $H_{(\alpha,a+\cc^\nu_\beta)}^q
H_{(b-\cc^\alpha_\mu,\beta)}^q$, $H_{(\alpha,a+\cc^\nu_\beta+1)}^q
H_{(b-\cc^\alpha_\mu-1, \beta)}^q$, $H_{(\alpha,b+\cc^\nu_\beta)}^q
H_{(a-\cc^\alpha_\mu,\beta)}^q$, or $H_{(\alpha,b+\cc^\nu_\beta-1)}^q
H_{(a-\cc^\alpha_\mu+1, \beta)}^q$. Let $H_{(\theta,c)}^q
H_{(d,\gamma)}^q$ be one of these terms after it has been rewritten
using Corollary \ref{co:Bott} so
that $(\theta,c)\in\DIW{k}$ and $(d,\gamma)\in\DIW{n}$.

Consider $H_{(\theta,c)}^q H_{(d,\gamma)}^q$ of the first form
since the others follow from essentially the same remark.  We verify
that $d-c < b-a$ (unless $H_{(\theta,c)}^q H_{(d,\gamma)}^q =
H_{(\mu,a)}^q H_{(a+1,\nu)}^q$) and hence by induction the theorem is
true since if
$d-c \leq 0$ then $(\theta,c,d,\gamma)\in\DIW{k+n}$.  By definition,
\begin{equation*}
\begin{split}
c &= \min_{1\le i\le k-1} \left\{a+\cc^\nu_\beta,\alpha_i+(k-i) \right\} \\
d &= \max_{1\le i\le n-1} \left\{b-\cc^\alpha_\mu, \beta_{i}-i\right\}
\end{split}
\end{equation*}
We note that
\begin{equation*}
c \geq \min_{1\le i\le k-1} \left\{a+\cc^\nu_\beta, \mu_i+(k-i) \right\}
 \geq \min \left\{a+\cc^\nu_\beta, a+1 \right\}
\end{equation*}
and hence $c \ge a$ with equality if and only if $\cc^\nu_\beta=0$.
Similarly
\begin{equation*}
d \le \max_{1\le i\le n-1} \left\{b-\cc^\alpha_\mu, \nu_{i}-i \right\}
\leq  \max \left\{b-\cc^\alpha_\mu, b-i \right\}
\end{equation*} 
and so $d \leq b$ with equality if and only if $\cc^\alpha_\mu=0$.
\end{proof}

The generalization of this statement, that $H_{\gamma^{(1)}}^q
H_{\gamma^{(2)}}^q \dotsm H_{\gamma^{(t)}}^q$ with $\gamma^{(i)} \in
\Z^{\eta_i}$ is in the
$\Z[q]$ span of the operators $H_{\alpha^{(1)}}^q
H_{\alpha^{(2)}}^q \dotsm H_{\alpha^{(t)}}^q$ with $\alpha^{(i)}\in
\Z^{\eta_i}$ and $(\alpha^{(1)},\alpha^{(2)}, \ldots, \alpha^{(t)})$
a dominant weight, is conjectured to be true \cite{KKSW}, but
does not seem to follow easily from these relations.

By setting $\ell=1$ and $m=0$ in \eqref{eq:comm}, 
we immediately obtain another relation among the operators
which gives us another basis theorem.

\begin{lem} \label{le:move}
For all $a\in\Z$, $\mu\in\DIW{k}$, and $\nu \in \DIW{n}$,
\begin{equation}\label{eq:move}
\sum_{\substack{\beta\in\DIW{n} \\ \nu/\beta\in\V}}
 (-q)^{\cc^\nu_\beta} H^q_{(\mu,a +\cc^\nu_\beta)}
H_\beta^q =
\sum_{\substack{\alpha\in\DIW{k} \\ \alpha/\mu\in\V}}
  (-q)^{\cc^\alpha_\mu} H_\alpha^q H^q_{(a - \cc^\alpha_\mu,\nu)}.
\end{equation}
\end{lem}
\begin{proof}
{}From equation \eqref{eq:comm},
\begin{equation}
H(U^k,v)H(Z^n) \Omega[-q v^{\mone} Z] = H(U^k) H(v, Z^n) 
\Omega[-q v U^{\mone}]
\end{equation}
The stated identity is obtained by
taking the coefficient of $u^\mu v^a z^\nu$.
\end{proof}

Before proving Theorem \ref{th:support} we give an example.
\begin{ex} We wish to write the composition of operators
$H_{(53)}^q H_{(2)}^q$ as a linear combination of operators
$H_\alpha^q H_\beta^q$ with $\alpha \in \Z^1$ and $\beta \in \Z^2$.
\begin{equation*}
H_{(53)}^q H_{(2)}^q - q H_{(54)}^q H_{(1)}^q 
= H_{(5)}^q H_{(32)}^q - q H_{(6)}^q H_{(21)}^q
\end{equation*}
\begin{equation*}
H_{(54)}^q H_{(1)}^q - q H_{(55)}^q H_{(0)}^q 
= H_{(5)}^q H_{(41)}^q - q H_{(6)}^q H_{(31)}^q
\end{equation*}
\begin{equation*}
H_{(55)}^q H_{(0)}^q = H_{(5)}^q H_{(50)}^q - q H_{(6)}^q H_{(40)}^q
\end{equation*}
Making repeated substitutions, we have the relation
\begin{align*}
H_{(53)}^q H_{(2)}^q = ~
&q^2 H_{(5)}^q H_{(50)}^q 
- q^3 H_{(6)}^q H_{(40)}^q  
+ q H_{(5)}^q H_{(41)}^q \\
&- q^2 H_{(6)}^q H_{(31)}^q 
+ H_{(5)}^q H_{(32)}^q 
- q H_{(6)}^q H_{(21)}^q
\end{align*}
\end{ex}

We now give the proof of Theorem \ref{th:support}.
\begin{proof} We prove that if $k>n$ then $\GGG(k,n)\subset \GGG(k-1,n+1)$
as the other part is proven in the same manner.
It is enough to show that if $k>n$,
$(\mu,a) \in \DIW{k}$ and $\nu\in\DIW{n}$ then
\begin{equation*}
H_{(\mu,a)}^q H_{\nu}^q \in \GGG_\ge(k-1,n+1)
\end{equation*}
We use \eqref{eq:move} with $k-1$ and $n$.
The proof proceeds by induction on $\mu_1 -a$.  If $\mu_1 - a = 0$,
then there is exactly one term on the left hand side of equation
\eqref{eq:move}, namely $H_{(\mu,a)}^q H_\nu^q$.  The right hand side
is a $\Z[q]$-linear combination of $H_\gamma^q H_\rho^q$
with $\gamma \in \Z^{k-1}$ and $\rho \in \Z^{n+1}$.

If $\mu_1 - a > 0$ then the left hand side is a linear combination of
operators of the form
$H_{(\gamma, c)}^q H_\beta^q$ with $\gamma_1 = \mu_1$ and 
\begin{align*}
c&=
\min_{1\le i\le k-1} \left\{ a+ \cc^\nu_\beta, \mu_i + (k-i)\right\} \\
&\ge \min_{1\le i\le k-1} \left\{ a+\cc^\nu_\beta, a + (k-i)\right\}.
\end{align*}
It follows that $c \ge a$ with equality if and only if $\nu = \beta$.
The right hand side of this equation is again an element of
$\GGG(k-1,n+1)$.
\end{proof}

Finally we prove Theorem \ref{th:comm}.
\begin{proof} We shall prove $\GGG(k+\ell,k)\subseteq\GGG(k,k+\ell)$
for $k,\ell>0$, the other inclusion being similar.
Consider \eqref{eq:comm} with $m=0$ and $n=k$ so $U=U^k$,
$V=V^\ell$, $Z=Z^k$:
\begin{equation} \label{eq:big move}
  H(U,V) H(Z) \Omega[-q V^\mone Z] =
  H(U) H(V,Z) \Omega[-q U^\mone V].
\end{equation}
Take the coefficient of $u^\alpha v^\beta z^\gamma$ in this equation,
where $(\alpha,\beta)\in\DIW{k+\ell}$ and $\gamma\in\DIW{k}$.
The entire right hand side consists of terms in $\GGG(k,k+\ell)$.
Expanding the expression $\Omega[-q V^* Z] = \prod (1-q z_i/v_j)$,
the left hand side is in the set
\begin{equation*}
  H^q_{\alpha,\beta} H^q_\gamma +
  \Z[q] \sum_{\beta',\gamma'}
  H^q_{\alpha,\beta'} H^q_{\gamma'}
\end{equation*}
where $\beta'\in\Z^\ell$ such that $\beta'_j \le \beta_j + k$ for
all $1\le j \le \ell$ and $|\beta'|>|\beta|$.  Rewriting a typical term
$H^q_{\alpha,\beta'}$ by Corollary \ref{co:Bott}, one obtains either $0$
or up to sign $H^q_{\alpha',\beta''}$, say.  This means
$\alpha'_1 + (k+\ell-1)$ is the
largest part of the weight $(\alpha,\beta')+\rho^{(k+\ell)}$.
Now $\alpha_1 \ge \beta_j$ for all $1\le j\le \ell$ 
by the dominance of $(\alpha,\beta)$, so
\begin{equation*}
\begin{split}
\alpha_1+k+\ell-1 &\ge \beta_j + k + \ell - 1 \ge \beta'_j +\ell-1 \\
&\ge \beta'_j + \ell-j.
\end{split}
\end{equation*}
It follows that $\alpha'_1=\alpha_1$, $(\alpha',\beta'')\in\DIW{k+\ell}$
and $|\alpha'|+|\beta''|> |\alpha|+|\beta|$.  There are only finitely
many elements of $\DIW{k+\ell}$ with these properties, so the terms
$H^q_{(\alpha',\beta'')} H^q_{\gamma'}$
can again be rewritten in the same way, and the process terminates.
\end{proof}

\section{Recovering identities via commutation relations}

The commutation relations in the previous section may be used
to recover some identities among the operators $H^q_\nu$
that correspond to known identities among generalized Kostka
polynomials.

\begin{prop} \label{pp:same width} For $a\in\Z$ and
positive integers $k$ and $n$,
\begin{equation} \label{eq:same width}
H^q_{(a^n)} H^q_{(a^k)} = H^q_{(a^k)} H^q_{(a^n)}
\end{equation}
\end{prop}
\begin{proof} In the proof of Theorem \ref{th:comm}, take
$\alpha=(a^k)$, $\beta=(a^\ell)$, and $\gamma=(a^k)$.  The
proof shows there is only one term on
either side of the identity coming from \eqref{eq:big move},
which agrees with \eqref{eq:same width} when $n=k+\ell$.
\end{proof}

\begin{prop} \label{re:one more}
\begin{equation} \label{eq:one more}
  H^q_{(a^k)} H^q_{((a+1)^k)} = q^k H^q_{((a+1)^k)} H^q_{(a^k)}.
\end{equation}
\end{prop}
\begin{proof} When Lemma \ref{le:move} is applied with
$\mu=(a^{k-1})$, $a$, and $\nu=((a+1)^k)$, there are no surviving terms on
the right and exactly two on the left, indexed by $\beta=\nu=((a+1)^k)$
and $\beta=(a^k)$.  Applying Corollary \ref{co:Bott} to the
term with $\beta=(a^k)$ the desired identity is obtained.
\end{proof}

\begin{prop} \label{pp:quad} For all $a\in\Z$ and $k\ge1$,
\begin{equation} \label{eq:quad}
  H^q_{(a^k)} H^q_{(a^k)} = H^q_{(a^{k+1})} H^q_{(a^{k-1})} +
        q^k H^q_{((a+1)^k)} H^q_{((a-1)^k)}.
\end{equation}
\end{prop}
\begin{proof} Apply Lemma \ref{le:move} with $n=k-1$, $\mu=(a^k)$,
and $\nu=(a^{k-1})$.  On the left side one has the single
nonvanishing term $H^q_{(a^{k+1})}H^q_{(a^{k-1})}$ corresponding to
the summand $\beta=\nu=(a^{k-1})$.  The right side has two nonvanishing
terms indexed by $\alpha=\mu=(a^k)$ and $\alpha=\mu+(1^k)=((a+1)^k)$.
The first of these terms is $H^q_{(a^k)} H^q_{(a^k)}$.
The second is $(-q)^k H^q_{((a+1)^k)} H^q_{(a-k,a^{k-1})}$.
But $H^q_{(a-k,a^{k-1})} = (-1)^{k-1} H^q_{((a-1)^k)}$ by
Corollary \ref{co:Bott}.
\end{proof}

This identity can be viewed as the trace of a short exact sequence that
resolves the ideal of a nilpotent conjugacy class closure over
the coordinate ring of a minimally larger one \cite{KKSW}.
The version of this identity for the fermionic form of 
generalized Kostka polynomials appears in \cite{ScWa}.
This $q$-character identity is put in a more general context
in \cite{HKOTY}.

\section{Generalization of Garsia-Procesi defining recurrence for
Kostka-Foulkes polynomials}

Manipulations of the definition allow us to derive commutation
relations between the $H_\nu^q$ and the operators $e_k^\perp$.
As a consequence we obtain a generalization of a defining recurrence
for the Kostka-Foulkes polynomials given in \cite{GP}.

Let $E(u)$ be the generating function of operators on $\La$ defined by
\begin{equation}
E(u) P[X] = P[X-u].
\end{equation}
By \eqref{eq:coprod perp} we have
\begin{equation*}
  P[X-u] = \sum_{\la\in\Part} s_\la^\perp(P)[X] s_\la[-u]
  = \sum_{k\ge0} e_k^\perp(P)[X] (-u)^k.
\end{equation*}
In other words
\begin{equation} \label{eq:e perp pleth}
  e_k^\perp P[X] = (-1)^k P[X-u]\coeff_{u^k}.
\end{equation}
The commutation relation of $H(Z^n)$ and $E(u)$ is:
\begin{equation*}
\begin{split}
  E(u) H(Z^n) P[X] &= E(u) P[X-(1-q)Z^\mone]\Omega[XZ]R(Z^n)  \\
  &= P[X-u-(1-q)Z^\mone]\Omega[(X-u)Z]R(Z^n) \\
  &= \Omega[-uZ] H(Z^n) E(u) P[X]
\end{split}
\end{equation*}
Taking the coefficient of $(-u)^k z^\la$ on both sides of this equation
we obtain the following relation.

\begin{prop} \label{pp:col H} Let $\la\in\DIW{n}$.  Then
\begin{equation*}
e_k^\perp H_\la^q = \sum_{\substack{\beta\in\DIW{n} \\ \la/\beta \in \V}}
  H_\beta^q e_{k - |\la| + |\beta|}^\perp
\end{equation*}
\end{prop}

Let $\eta=(\eta_1,\eta_2,\dotsc,\eta_t)$ be a fixed sequence
of positive integers summing to $n$.  For any weight $\gamma\in\Z^n$,
write $\gamma^{(1)}\in\Z^{\eta_1}$ for the first $\eta_1$ parts of
$\gamma$, $\gamma^{(2)}\in\Z^{\eta_2}$ for the next $\eta_2$ parts of
$\gamma$, etc.

\begin{prop} \label{pp:col skew} Let $k$ be a fixed positive integer,
$\eta=(\eta_1,\eta_2,\dotsc,\eta_t)$ a sequence of positive
integers summing to $n$, $\alpha\in\DIW{n}$, and $\gamma\in\Z^n$ such that
$\gamma^{(i)}\in\DIW{\eta_i}$ for all $i$ and $|\gamma|-|\alpha|=k$.  Then
\begin{equation*}
  \sum_{\substack{\nu\in\Z^n \\ |\gamma|-|\nu|=k \\
    \nu^{(i)}\in\DIW{\eta_i} \\ \gamma^{(i)}/\nu^{(i)} \in \V }}
        K_{\alpha,\nu,\eta}(q) =
  \sum_{\substack{\la\in\DIW{n} \\ |\la|-|\alpha|=k \\ \la/\alpha\in\V }}
        K_{\la,\gamma,\eta}(q)
\end{equation*}
\end{prop}
\begin{proof} Let $\eta$ and $\gamma$ be as above.  For the moment assume
that the entries of $\gamma$ are positive.
Apply $e_k^\perp$ to \eqref{eq:genfun} and
apply Proposition \ref{pp:col H} to both sides to commute $e_k^\perp$
to the right of the $H$ operators: 
\begin{equation*}
  \sum_{\substack{\nu\in\Z^n \\ \nu^{(i)} \in\DIW{\eta_i} \\
  \gamma^{(i)}/\nu^{(i)} \in\V }}
  H_{\nu^{(1)}}^q H_{\nu^{(2)}}^q \cdots H_{\nu^{(t)}}^q
        e_{k -|\gamma|+|\nu|}^\perp = 
  \sum_{\la\in\DIW{n}}
  \sum_{\substack{\beta\in\DIW{n} \\ \la/\beta\in\V}}
        K_{\la,\nu,\eta}(q) H_\beta^q e_{k-|\la|+|\beta|}^\perp.
\end{equation*}
Since $\gamma$ has positive parts, all the $\nu$ have nonnegative parts.
Expanding the left hand side using \eqref{eq:genfun} and applying the
resulting operators to $1\in\La$, we obtain
\begin{equation*}
\sum_{\substack{\nu\in\Z^n \\ |\gamma|-|\nu| = k \\
        \nu^{(i)} \in\DIW{\eta_i} \\ \gamma^{(i)}/\nu^{(i)} \in\V}}
        \sum_{\mu\in\DIW{n}} K_{\mu,\nu,\eta}(q) H^q_\mu \, 1 =
\sum_{\la\in\DIW{n}} \sum_{\substack{\beta\in\DIW{n} \\
        |\la|-|\beta|=k \\ \la/\beta\in\V}}
  K_{\la,\nu,\eta}(q) H_\beta^q \, 1.
\end{equation*}
Assume for the moment that $\alpha$ has nonnegative parts.
For such $\alpha$, $H^q_\alpha \, 1 = s_\alpha[X]$.
Taking the coefficient of $s_\alpha[X]$ on both sides,
we obtain the desired relation.

The statement is true in general since
$K_{\la-(a^n),\gamma-(a^n),\eta}(q) = K_{\la,\gamma,\eta}(q)$
for all integers $a$.
\end{proof}

The recurrence for the Kostka-Foulkes polynomials given in
\cite{GP} is recovered by setting $\eta=(1^n)$
in Proposition \ref{pp:col skew},
as $H^q_{(m)}$ in our notation is $H_m$ in theirs.
In the case $\eta=(1^n)$ the situation is particularly nice;
the nondominant operators can be made dominant using the relation
$H^q_{(m)} H^q_{(m+1)} = q H^q_{(m+1)} H^q_{(m)}$.
For general $\eta$ the identities needed to rewrite the
nondominant terms as dominant ones, can get complicated
and produce negative signs.

\section{Jing's operators}

We give corresponding constructions for Jing's Hall-Littlewood
vertex operators \cite{J} \cite[Ex. III.5.8]{M}.  Define
\begin{equation*}
\begin{split}
  \BB(Z^k) P[X] &= P[X - Z^\mone] \Omega[X Z (1-q)]  \\
  B(Z^k) P[X] &= P[X- Z^\mone] \Omega[X Z (1-q)] R(Z^k) \\ 
  B^q_\nu &= B(Z^k)\coeff_{z^\nu}  
\end{split}
\end{equation*}
for $\nu\in\Z^k$.  The generating series of operators
$B(z)$ (for a single variable $z$) and $\BB(Z^k)$ defined
here, coincide with the operators $B(z)$ and $B(z_1,z_2,\dots,z_k)$
in the notation of \cite[Ex. III.5.8]{M}.
$B(z)$ in our notation is the operator $H(z)$ in Jing's \cite{J}.
However the operators $B_\nu$ themselves are not studied
by Jing or Macdonald.

Let $F$ be the plethystic operator $F P[X] = P[X(1-q)]$, with inverse
operator $F^{-1} P[X] = P[X/(1-q)]$.  Then it is not hard to show that
\begin{equation}\label{B H}
  B(Z^k) = F \circ H(Z^k) \circ F^{-1}.
\end{equation}

This gives the following analogue of \eqref{eq:genfun} for $B(Z^k)$.
\begin{equation}  \label{eq:genfun B}
B_{\gamma^{(1)}}^q B_{\gamma^{(2)}}^q \cdots B_{\gamma^{(t)}}^q =
 \sum_{\la\in\DIW{n}} K_{\la,\gamma,\eta}(q) B_\la^q.
\end{equation}

\end{document}